\documentclass{amsart}
\usepackage{amssymb,latexsym,amsmath}

\theoremstyle{plain}
\newtheorem{lemma}{Lemma}[section]
\newtheorem{theorem}[lemma]{Theorem}
\newtheorem{corollary}[lemma]{Corollary}
\newtheorem*{main}{Main theorem}
\theoremstyle{definition}
\newtheorem{definition}[lemma]{Definition}
\theoremstyle{remark}

\begin{document}
\title[]
    {The relative hyperbolicity of one-relator relative
    presentations}
\author[]{Le Thi Giang}
\address{Faculty of mechanics and mathematics\\
        Moscow State University\\
        Moscow 119992, Leninskie gory, MSU}
\email{giangmsu@gmail.com}
\begin{abstract}
We prove that if $G$ is a free-torsion group and $w(t)$ is a word
in the alphabet $G \sqcup \{t^{\pm 1}\}$ with exponent sum one,
then the group $\langle G,t|(w(t))^k = 1 \rangle$, where $k \geq
2$, is relatively hyperbolic with respect to $G$.
\end{abstract}

\maketitle

\section{Introduction}
Let us recall that a group $\bar{G}$ is called relatively
hyperbolic with respect to a collection of its subgroups
$\{H_\lambda\}_{\lambda \in \Lambda}$ if $\bar{G}$ admits a
relatively finite presentation
$$\langle X,\{H_\lambda\}_{\lambda \in \Lambda} | R = 1 , R \in \mathcal{R}\rangle ,$$
where $X$ and $\mathcal{R}$ are finite, satisfying a linear
relative isoperimetric inequality. That is, there exists $C > 0$
such that for every word $w$ in the alphabet $X^{\pm 1} \cup
\mathcal{H}$, where $\mathcal{H} = \sqcup_{\lambda \in
\Lambda}(H_{\lambda} \backslash \{1\})$, representing the identity
in the group $\bar{G}$, there exists an expression
$$w =_F \prod_{i=1}^k f_i^{-1}R_i^{\pm 1}f_i,$$
with the equality in the group $F =(*_{\lambda \in
\Lambda}H_{\lambda})*F(X)$, where $F(X)$ is the free group with
basis $X$,$R_i \in \mathcal{R}, f_i \in F$ and $k \leq C\|w\|$,
where $\|w\|$ is the length of the word $w$.

The following theorem belongs to B. Newman~\cite{Newman68}, see
also~\cite{LS77}.

\begin{theorem}
    Any one-relator group
    $$ \langle a,b,\ldots |R^k \rangle, $$
    where $k > 1$, is hyperbolic.
\end{theorem}

We obtain a similar fact for one-relator relative presentations.

Let $G$ be a torsion-free group and  let $t$ be a letter, $t
\notin G$.

\begin{definition}
    The word $w = \prod_{i=1}^{n}
    g_it^{\epsilon_i} \in G*{\langle t\rangle}_\infty$, where
    $\epsilon_i = \pm 1$, is called unimodular if $\sum \epsilon_i =
    1$.
\end{definition}

\begin{main}\label{T:mainthe}
    A group
    $$\tilde{G}=\langle G,t| w^k =1\rangle  := G*{\langle t\rangle}_{\infty}/{\ll}w^k{\gg},$$
    where the word $w$ is unimodular, $k \geq 2$ and $G$ is a
    torsion-free group, is relatively hyperbolic with respect to $G$.
\end{main}

The following example shows that if $w$ is not unimodular and $G$
is not torsion-free group then the theorem is not true.

\textbf{Example.} Suppose that $G = \langle a \rangle_3$ and
$\tilde{G} = \langle G,t| (a^{-1}t^{-1}at)^3=1 \rangle$. Then the
elements $u = t^{-1}ata$, $v = at^{-1}at$ generate an abelian free
group of rank 2, so the group $\tilde{G}$ can not be relatively
hyperbolic with respect to $G$ (see \cite{Osin06}).

Recall that a group $H$ is called \emph{SQ-universal} if every
countable group can embedded into a quotient of $H$. In
\cite{AMO07} it was proved the SQ-universality of non-elementary
properly relatively hyperbolic groups.

\begin{corollary}\label{C:SQuniv}
    Under the conditions of theorem \ref{T:mainthe} the group
    $\tilde{G}$ is SQ-universal.
\end{corollary}

Theorem \ref{T:mainthe} implies the following corollary (see
\cite{Osin06}).

\begin{corollary}\label{C:WordProblem}
    If the word problem is solvable in the group $G$, then it is solvable in
    $\tilde{G}$.
\end{corollary}

One-relator groups have been well studied. In these groups the
word problem is solvable \cite{LS77},\cite{MKS66}. In \cite{How91}
J. Howie studied the quotient of the free products of the type $G
= (A*B)/r^n,$ where $n > 3$. For this type of groups the word
problem is also solvable. In \cite{Ar06} A. Juhasz generalized
this result for some amalgamated free products. These results were
obtained by using small cancellation theory and it is unknown
whether these results can be extended to the cases $n =2,3$.

In this work we use Howie diagrams over a relative presentation
and  Klyachko's method of multiple motion on diagrams.
\\
\textbf{Agreement.} Thought out this paper, the symbols
$G,\tilde{G},w,k,t$ have the same sense as above.
\\
\textbf{Acknowledgements.} I am grateful to my advisor Ant. A.
Klyachko, who helped me so much. I also thank prof. A. Yu.
Olshanskii for many useful comments.

\section{Algebraic lemmata}\label{S:Lemmata}

\begin{lemma}\label{L:Inj}
    The natural mapping $G\to \tilde{G}$ is injective.
\end{lemma}

\emph{Proof.} In \cite{Kl93} it was proved that, if a word $w$ is
unimodular, then the natural mapping $G \to \hat{G}=\langle G,t |
w =1\rangle := G*{\langle t\rangle}_{\infty}/ {\ll}w{\gg} $ is
injective. Now, consider the natural mappings $G \to \tilde{G} \to
\hat{G}$. If $G$ does not embed into $\tilde{G}$ then $G$ does not
embed into $\hat{G}$. This contradiction proves Lemma~\ref{L:Inj}.

The following lemma is similar to Lemma 2 from \cite{Kl05}.

\begin{lemma}\label{L:Otpre}\cite{Kl05}
    Let the word $w$ be unimodular and nonconjugate to a word  of the
    form $gt$, where $g \in G$.
    Then the group $\tilde{G}$ has a (relative) presentation of the form
    \begin{equation}\label{E:equ1}
    \tilde{G} = \langle H,t| \{p^t=p^{\varphi},p \in P\backslash
    \{1\}\},(ct\prod_{i=0}^{m}(b_i{a_i}^t))^k = 1\rangle,
    \end{equation}
    where $a_i, b_i, c \in H, P$ and
    $P^{\varphi}$ are isomorphic subgroups of $H$, $\varphi : P \to
    P^{\varphi}$ is an isomorphism  and the following conditions hold
    \begin{enumerate}
    \item[1)] $m \geq 0$ (i.e, the product in \textup{(\ref{E:equ1})} is nonempty);
    \item[2)] $a_i \notin P,b_i \notin P^{\varphi} ;$
    \item[3)] $\langle P,a_i\rangle = P*{\langle a_i\rangle}_{\infty},\langle P^{\varphi},
    b_i\rangle = P^{\varphi}*{\langle b_i\rangle}_{\infty}$ in
    $H$;
    \item[4)] The groups $H,P,P^{\varphi}$ are free products of finitely
    many isomorphic copies of $G$: $H = G^{(0)}*...*G^{(s)}$, $P =
    G^{(0)}*...*G^{(s-1)}$, $P^{\varphi}= G^{(1)}*...*G^{(s)}$, where $s
    \geq 0$ (for $s = 0$,  $P$ and $P^{\varphi}$ are assumed to be
    trivial), and the isomorphism $\varphi$ is the shift:
    $(G^{(i)})^{\varphi} = G^{(i+1)}$.
    \end{enumerate}
\end{lemma}

\section{Howie diagrams}\label{S:Diag}
Consider a map on an oriented two-dimensional sphere. The corners
of this map are labelled by elements of some group $H$. The edges
are directed and labelled by the letter $t$.

The \emph{label of a vertex} is defined as the product of the
labels of all corners near this vertex listed clockwise. The label
of a vertex is an element of $H$ defined up to conjugation.

To obtain the \emph{label of a face}, we should walk along its
boundary anticlockwise and write down the labels of all its
corners and edges; the label of an edge should be written as its
inverse if we walk through it against the arrow. The label of a
face is an element of the group $H*\langle t \rangle_\infty$ (the
free group of $H$ and the cyclic group with basis $t$), defined up
to a cyclic permutation.

Such a labelled map is called a \emph{Howie disk diagram} over a
relative presentation
\begin{equation}\label{E:equ2}
\langle H,t | w_1=1, w_2=1, \ldots \rangle,
\end{equation}
if
\begin{itemize}
\item one face is separated out and called \emph{exterior}, the
remaining faces are called \emph{interior};
\item the label of each interior face is a cyclic permutation of
one of the words ${w_i}^{\pm1}$;
\item the label of each vertex is the identity element of $H$.
\end{itemize}
The diagram is said to be \emph{reduced} if it contains no such
edge $e$ that both faces containing $e$ are interior, these faces
are different and their labels written starting with the label of
edge $e$ are mutually inverse (such a pair of faces with a common
edge is called a \emph{reducible pair}).

The following lemmata are the key ingredients of the proof in the
rest of the paper.

\begin{lemma}\label{L:How83}~\textup{\cite{How83}}
    If the natural mapping from the group $H$ to the group with relative
    presentation~\textup{(\ref{E:equ2})} is injective, then the image of
    an element $u \in H*\langle t\rangle_\infty \backslash \{1\}$ is
    identity in the group~\textup{(\ref{E:equ2})} if and only if there
    exists a disk diagram over this presentation whose exterior face
    is labelled by $u$. A minimal (with respect to the
    number of faces) such diagram is reduced.
\end{lemma}

Let $\varphi : P \to P^{\varphi}$ be an isomorphism between two
subgroups of a group $H$. A relative presentation of the form
\begin{equation}\label{E:equ3}
\langle H,t | \{ p^t = p^{\varphi};p \in P\backslash \{1\}\},w_1 =
1, w_2 = 1,...\rangle,
\end{equation}
is called a $\varphi-$\emph{presentation}. A diagram over a
$\varphi-$presentation (3) is called $\varphi-$\emph{reduced} if
it is reduced and different interior cells with labels of the form
$p^tp^{-\varphi}, p \in P$, have no common edges.

\begin{lemma}\label{L:Kl051}~\textup{\cite{Kl05}}
    The complete $\varphi$-analog of lemma~\textup{\ref{L:How83}}
    is valid, i.e any minimal diagram over
    presentation~\ref{E:equ3} whose exterior face is labelled by
    $u$ is $\varphi-$reduced.
\end{lemma}

\emph{A multiple motion of period} $T \in \mathbb{R}$ on a diagram
with faces $\{D_i, i \in I\}$ is called a set of mappings $\{
\alpha_{i,j}:\mathbb{R} \to
\partial D_i;i \in I, j = 1,...,d_i \}$ (called cars),
satisfying the following periodicity conditions:
\begin{enumerate}
\item[1)] $\alpha_{i,j}(t+T) = \alpha_{i,j+1}(t)$ for any $t \in
\mathbb{R}$ and $j \in \{ 1,...,d_i \}$ (subscripts are modulo
$d_i$);
\item[2)] there exists such a partitioning of each circle ${\partial}D_i$ into $d_i$
arcs that during the time interval $[0,T]$ each car $\alpha_{i,j}$
moves along the  $j$th arc.
\end{enumerate}
The positive integers $d_i$ are called the \emph{multiplicities}
of the multiple motion.

The multiple motion $\alpha_{i,j}$ is called a \emph{multiple
motion with separated stops} if every car moves without U-turns
and infinite decelerations and accelerations moving around the
boundary of its face anticlockwise, possible stopping for a finite
time at some corners; and there exists a set of corners called the
stop corners such that:
\begin{enumerate}
\item[1)] the cars stop only at stop corners (possibly, at some stop
corners, the cars do not stop);
\item[2)] at each vertex $v$ having stop corners at it, the stops are
separated in the following sense : let $c_1,...,c_k$ be all stop
corners at $v$ enumerated anticlockwise; it is required that, for
each $i$, at corners  $c_i$ and $c_{i+1}$ (subscripts are modulo
$k$) cars are never located simultaneously (in particular, this
implies that $k \geq 2$).
\end{enumerate}

If the number of cars being at a moment $t \in \mathbb{R}$ at a
point $p$ of the sphere equals the multiplicity of this point (in
other words, either two cars simultaneously pass the same internal
point $p$ of an edge at the moment $t$ or there are $k$ cars at a
vertex $p$ of degree $k$ at the moment $t$), then we say that a
\emph{complete collision} occurs at the point $p$ at the moment
$t$; the point $p$ is called a \emph{point of complete collision}.

\begin{lemma}\label{L:Kl052}\textup{\cite{Kl05}}
    The number of points of complete collision of a multiple motion
    with separated stops on a sphere cannot be smaller than
    $$2 + \sum_{i \in I}(d_i -1),$$
    where $d_i$ are the multiplicities of the multiple motion.
\end{lemma}

\section{Standard multiple motion}\label{S:Ctand}
In this section we will define a multiple motion on Howie diagram
over the presentation~(\ref{E:equ1}). \\
A Howie diagram can have corners of four kinds: $(+ +), (- -), (+
-)$ and $(- +)$ (fig.2).\\

\begin{picture}(200,50)
\put(0,0){\vector(1,0){20}} \put(10,0){\line(1,0){30}}
\put(0,50){\vector(0,-1){30}} \put(0,0){\line(0,1){20}}
\put(0,0){\circle*{5}} \put(15,20){$(++)$}
\put(70,0){\line(1,0){20}} \put(120,0){\vector(-1,0){30}}
\put(70,0){\vector(0,1){30}}
\put(70,30){\line(0,1){20}}\put(70,0){\circle*{5}} \put(90,20){$(-
-)$} \put(190,0){\vector(-1,0){30}} \put(140,0){\line(1,0){20}}
\put(140,50){\vector(0,-1){30}} \put(140,0){\line(0,1){20}}
\put(140,0){\circle*{5}} \put(160,20){$(+ -)$}
\put(210,0){\vector(1,0){30}} \put(240,0){\line(1,0){20}}
\put(210,0){\vector(0,1){30}} \put(210,30){\line(0,1){20}}
\put(210,0){\circle*{5}} \put(230,20){$(- +)$} \put(300,0){fig. 2}
\end{picture}\\

A vertex of the kind

\begin{picture}(200,50)
\put(150,20){\circle*{5}} \put(120,20){\vector(1,0){20}}
\put(120,20){\line(1,0){30}} \put(180,20){\vector(-1,0){20}}
\put(150,20){\line(1,0){30}} \put(130,0){\vector(1,1){15}}
\put(130,0){\line(1,1){20}}  \put(170,0){\vector(-1,1){15}}
\put(150,20){\line(1,-1){20}} \put(150,50){\vector(0,-1){20}}
\put(150,20){\line(0,1){30}} \put(300,20){fig. 3}
\end{picture}\\
is called \emph{a sink},\\
and a vertex of the kind

\begin{picture}(200,60)
\put(150,30){\circle*{5}} \put(150,30){\vector(0,1){30}}
\put(150,30){\vector(1,0){30}} \put(150,30){\vector(-1,-1){30}}
\put(150,30){\vector(-1,0){30}} \put(150,30){\vector(1,-1){30}}
\put(300,30){fig. 4}
\end{picture}\\
is called \emph{a source}.\\
\\
The following lemma is obvious.

\begin{lemma}\label{L:l51}
    In the anticlockwise listing of the corners at a
    vertex $v$, the corners of type $(++)$ alternate with
    corners of type $(- -)$. If at a vertex $v$ there are no
    corners of type $(++)$, or, equivalently, there are no
    corners of type $(- -)$ then either all corners at $v$ are
    of type $(+ -)$ (the vertex is a sink),or all corners at $v$ are of
    type $(- +)$(the vertex is a source).
\end{lemma}

Now we define a standard multiple motion on the interior cells of
a disk diagram over presentation~(\ref{E:equ1}), which is similar
to the standard motion in \cite{Kl05}.
\begin{itemize}
\item the car going around a face with label $p^{-\varphi}p^t$ moves anticlockwise
uniformly with unit speed (one edge per unit time) visiting the
corner of type $(+ -)$ at the even moments of time (fig.5a);
\item on the boundary of a face with label $(ct\prod_{i=0}^mb_ia_i^t)^{\pm
k}$ there are $k$ moving cars. Each of them at moment zero is in a
corner with label  $b_0^{\pm 1}$ and moves anticlockwise during
the time intervals $[0,4m+2]$ on its own arc with label
$(ct\prod_{i=0}^mb_ia_i^t)^{\pm 1}$;
\item for $m > 0$ the car moving on an arc with label
$ct\prod_{i=0}^mb_ia_i^t$ stays at the corner of type $(++)$
during the time intervals $[2m + 2, 4m + 1] + (4m + 2)\mathbb{Z}$
and moves uniformly with unit speed all the remaining time; for $m
= 0$ such a car moves without stops with speed 2 on the positive
edges and with speed 1 on the negative ones (fig.5b);
\item for $m > 0$ the car moving on an arc with label $(ct\prod_{i=0}^mb_ia_i^t)^{-1}$
stays at the corners of type $(- -)$ during the time intervals
$[1, 2m] + (4m + 2)\mathbb{Z}$ and moves uniformly with unit speed
all the remaining time; for $m = 0$ such a car moves without stops
with speed 2 on the negative edges and with speed 1 on the
positive ones (fig.5c).
\end{itemize}

\newpage
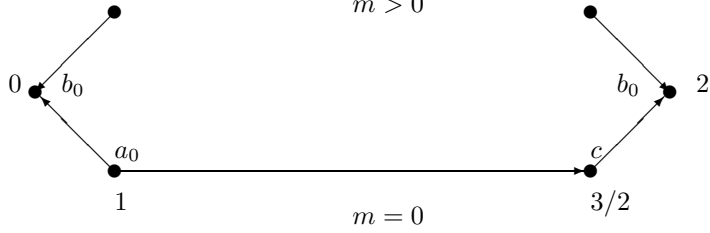
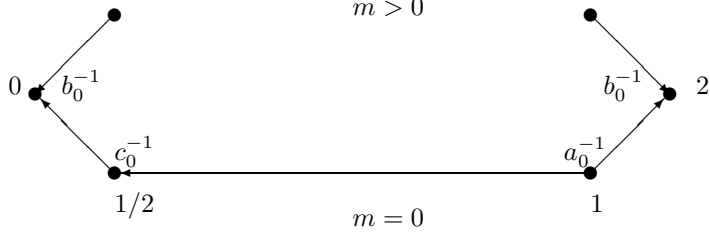
\begin{picture}(300,50)
\put(100,20){\circle*{5}} \put(250,20){\circle*{5}}
\put(175,20){\oval(150,50)} \put(110,20){$p^{-\varphi}$}
\put(175,45){\vector(-1,0){3}} \put(240,20){$p$}
\put(175,-5){\vector(-1,0){3}} \put(50,20){$0,2,4,...$}
\put(260,20){$1,3,5,...$} \put(300,20){fig.5a}
\end{picture}

\begin{picture}(300,200)
\put(10,150){\circle*{5}} \put(25,120){\circle*{5}}
\put(25,120){\vector(-1,2){13}} \put(25,180){\circle*{5}}
\put(25,180){\vector(-1,-2){13}} \put(0,145){$0$}
\put(20,145){$b_0$} \put(25,100){$1$} \put(25,125){$a_0$}
\put(60,120){\circle*{5}} \put(25,120){\vector(1,0){33}}
\put(60,100){$2$} \put(60,125){$b_1$}
\put(95,120){\vector(-1,0){33}} \put(95,118){ - - - }
\put(150,120){\circle*{5}} \put(150,120){\vector(-1,0){33}}
\put(140,100){$2m-1$} \put(140,125){$a_{m-1}$}
\put(185,120){\circle*{5}} \put(150,120){\vector(1,0){33}}
\put(180,100){$2m$} \put(175,125){$b_m$}
\put(220,120){\circle*{5}} \put(220,120){\vector(-1,0){33}}
\put(210,100){$2m+1$} \put(210,125){$a_m$}
\put(220,120){\vector(1,0){33}} \put(255,120){\circle*{5}}
\put(255,100){$[2m+2,4m+1]$} \put(245,125){$c$}
\put(255,120){\vector(1,2){13}} \put(270,150){\circle*{5}}
\put(280,145){$4m+2$} \put(255,145){$b_0$}
\put(255,180){\circle*{5}} \put(255,180){\vector(1,-2){13}}
\put(150,80){$m > 0$} \put(150,0){$m = 0$}
\put(30,50){\circle*{5}} \put(20,50){$0$} \put(40,50){$b_0$}
\put(60,20){\circle*{5}} \put(60,20){\vector(-1,1){28}}
\put(60,5){$1$} \put(60,25){$a_0$} \put(60,20){\vector(1,0){178}}
\put(240,20){\circle*{5}} \put(240,5){$3/2$} \put(240,25){$c$}
\put(270,50){\circle*{5}} \put(240,20){\vector(1,1){28}}
\put(250,50){$b_0$} \put(280,50){$2$} \put(60,80){\circle*{5}}
\put(60,80){\vector(-1,-1){30}} \put(240,80){\circle*{5}}
\put(240,80){\vector(1,-1){30}} \put(300,70){fig.5b}
\end{picture}

\begin{picture}(300,200)
\put(10,150){\circle*{5}} \put(25,120){\circle*{5}}
\put(25,120){\vector(-1,2){13}} \put(25,180){\circle*{5}}
\put(25,180){\vector(-1,-2){13}} \put(0,145){$0$}
\put(20,145){$b_0^{-1}$} \put(15,100){$[1,2m]$}
\put(25,125){$c^{-1}$} \put(60,120){\circle*{5}}
\put(60,120){\vector(-1,0){33}} \put(55,100){$2m+1$}
\put(60,125){$a_m^{-1}$} \put(60,120){\vector(1,0){33}}
 \put(95,120){\circle*{5}}
\put(130,120){\vector(-1,0){33}} \put(95,100){$2m+2$}
\put(95,125){$b_m^{-1}$} \put(130,117){- - - - - - - -}
\put(220,120){\circle*{5}} \put(185,120){\vector(1,0){33}}
\put(210,100){$4m$} \put(210,125){$b_1^{-1}$}
\put(255,120){\vector(-1,0){33}} \put(255,120){\circle*{5}}
\put(255,100){$4m+1$} \put(245,125){$a_0^{-1}$}
\put(255,120){\vector(1,2){13}} \put(270,150){\circle*{5}}
\put(280,145){$4m+2$} \put(250,145){$b_0^{-1}$}
\put(255,180){\circle*{5}} \put(255,180){\vector(1,-2){13}}
\put(150,80){$m > 0$} \put(150,0){$m = 0$}
\put(30,50){\circle*{5}} \put(20,50){$0$} \put(40,50){$b_0^{-1}$}
\put(60,20){\circle*{5}} \put(60,20){\vector(-1,1){28}}
\put(60,5){$1/2$} \put(60,25){$c_0^{-1}$}
\put(240,20){\vector(-1,0){178}} \put(240,20){\circle*{5}}
\put(240,5){$1$} \put(230,25){$a_0^{-1}$}
\put(270,50){\circle*{5}} \put(240,20){\vector(1,1){28}}
\put(245,50){$b_0^{-1}$} \put(280,50){$2$}
\put(60,80){\circle*{5}} \put(60,80){\vector(-1,-1){30}}
\put(240,80){\circle*{5}} \put(240,80){\vector(1,-1){30}}
\put(300,70){fig.5c}
\end{picture}
\newpage
Let $u$ be the label of the exterior face. For $m > 0$  we can
suppose that on the exterior face there are at least one corner of
type $(- -)$ and one of type $(++)$ (it means that in the word $u$
there are  subwords of form $t^{-1}g_1t^{-1}$ and $tg_2t$, where
$g_1, g_2 \in G$, such subwords will always be found, if
necessary, conjugating $u$ by $t^n, n > 0$). There is one car
moving on the exterior face. During the time intervals $[1,2m] +
(4m + 2)\mathbb{Z}$ and $[2m+2,4m+1] + (4m + 2)\mathbb{Z}$ the car
stays at the corners of type $(- -)$ and $(++)$. This car moves
uniformly anticlockwise with the same speed all the remaining
time.

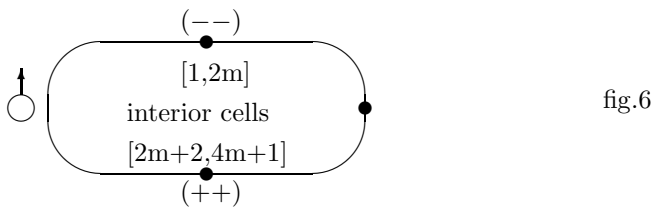
\begin{picture}(300,100)
\put(150,50){\oval(120,50)} \put(120,45){interior cells}
\put(300,50){fig.6} \put(80,50){\circle{10}}
\put(80,55){\vector(0,1){10}} \put(150,75){\circle*{5}}
\put(140,60){[1,2m]} \put(140,80){$(- -)$} \put(140,15){$(++)$}
\put(120,30){[2m+2,4m+1]} \put(150,25){\circle*{5}}
\put(210,50){\circle*{5}}
\end{picture}

For $m = 0$ the car moves uniformly with the same speed of period
2.

\begin{lemma}\label{L:Collision}
    The standard multiple motion on a diagram over relative
    representation \textup{(\ref{E:equ1})} is a motion with separated stops.
    In the interior cells the complete collisions can occur only at vertices
    being sources or sinks.
\end{lemma}
\emph{Proof.} Let all corners of types $(++)$ and $(- -)$ be the
stop corners. Then, the stops are separated by Lemma~\ref{L:l51}
and the fact that the schedule of the standard motion is such that
cars are never located at the same time at corners of type $(++)$
and $(- -)$. \\
A collision on an edge at a moment $t$ means that at this moment
the direction of the motion of one of the cars coincides with the
direction of the edge, while the direction of the motion of the
other colliding car is opposite to the direction of the edge
(fig.7)

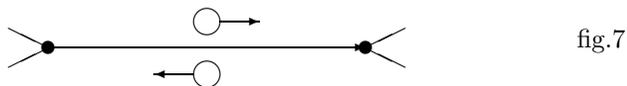
\begin{picture}(200,50)
\put(100,25){\vector(1,0){120}} \put(100,25){\circle*{5}}
\put(220,25){\circle*{5}} \put(100,25){\line(-2,1){15}}
\put(100,25){\line(-2,-1){15}}  \put(220,25){\line(2,-1){15}}
\put(220,25){\line(2,1){15}} \put(160,15){\circle{10}}
\put(155,15){\vector(-1,0){15}} \put(160,35){\circle{10}}
\put(165,35){\vector(1,0){15}} \put(300,25){fig.7}
\end{picture}

But the schedule of the standard multiple motion on interior cells
is such that, at each moment $t$, either all cars being on edges
move in the direction of the edge (this is so when the integer
part of $t$ is odd), or all cars being on edges move in the
direction opposite to the direction of the edge (this is so when
the integer part of $t$ is even). Therefore, collisions can occur
only at vertices; the separateness of stops implies that a vertex
of complete collision can not have stop corners and, therefore, is
a source or a sink. The lemma is proven.

\section{Proof of main theorem}\label{S:Proof}
By Lemma (\ref{L:Inj}), the natural mapping $G \to \tilde{G}$ is
injective. So for every word $u \in G*\langle t\rangle_{\infty} $,
where $u =_{\tilde{G}}1$, there is a reduced disk diagram over the
presentation of $\tilde{G}$ with a exterior face with label $u$
(Lemma \ref{L:How83}). We denote by $|u|$ the number of edges of
exterior face of this diagram.

It is well-known that in this case $u$ has a presentation
$$u = \prod_{i=1}^h {f_i}^{-1}w^{\pm k}f_i, $$
where h is the number of faces of the disk diagram.

We will prove the following inequality
\begin{equation}\label{E:equa4}
h < C|u|,
\end{equation}
where $C$ is a constant. In fact, this inequality is equivalent to
the inequality in the definition of hyperbolicity, because the
number of edges of the exterior face of the disk diagram exactly
equals to a half of the length of $u$.

We consider two cases:

\textbf{Case 1.}  $|w| > 1$. By Lemma~\ref{L:Otpre}, $\tilde{G}$
has a relative presentation~(\ref{E:equ1}). The construction of
$H$ implies that the natural mapping $H \to \tilde{G}$ is
injective, so for all words $u$ in the alphabet $H \cup \{t^{\pm
1}\}$, where $u =_{\tilde{G}}1$, there is a $\varphi-$reduced disk
diagram $M$ over the presentation~(\ref{E:equ1}) whose  exterior
face is labelled by $u$ (Lemma~\ref{L:Kl051}).

We define a standard multiple motion on $M$ as in
section~\ref{S:Ctand}. Let us show that all complete collisions of
the standard motion on this map occur only at the boundary of
exterior face.

Indeed, according to Lemma \ref{L:Collision}, a complete collision
on interior cells can occur only at vertices being sinks or
sources.

Suppose that a vertex of complete collision is a sink. Then, all
of corners at this vertex are of type $(+ -)$; the label of each
of these corners is either $p^\varphi$, where $p \in P$, or
${b_i}^{\pm 1}$. If at this vertex there are a corner labelled by
${b_i}^{\pm 1}$ and a corner labelled by ${b_j}^{\pm 1}$, where $i
\not= j$, then a complete collision does not occur at this vertex,
because these two corners are never visited at the same time.  If
at this vertex there is a corner labelled by ${b_i}^{\pm 1}$ and a
corner labelled by $p^\varphi$ then the label of a vertex being a
sink at which a complete collision occurs is
$$\prod_{j}({b_i}^{\epsilon_j}{p_j}^\varphi),$$
where $\epsilon_j \in \mathbb{Z}, p_j \in P$. But the label of an
vertex must be 1, thus we have a nontrivial (because the diagram
is $\varphi-$reduced) relation of the specified form, which
contradicts property 3 from Lemma~\ref{L:Otpre}.

For the case, when a vertex of complete collision is a source, we
can prove in the same way.

Suppose that $d$ is the number of interior faces with label of
type $p^tp^{-\varphi}$ (fig.5a), $e$ is the number of other faces
(fig.5b, fig.5c), $f$ is the number of edges of exterior face of
the disk diagram $M$.

By Lemma \ref{L:Kl052}, the number of complete collisions is not
smaller than $2 + \sum_{i = 1}^e (k - 1) = 2 + e(k-1)$. Besides
complete collision can occur only at the boundary of the exterior
face, therefore the number of complete collisions is not greater
than the number of edges of the exterior face of $M$ (because at
the boundary of the exterior face the car almost is located at
stop corners and in the remaining time moves very quickly, so at
each edge of the exterior face there is at most one point of
complete collision). Thus, we have $2+ e(k-1) \leq f$, i.e., $e
\leq \frac{f-2}{k-1} <\frac{f}{k-1}$.

So over the presentation~(\ref{E:equ1}) the word $u$ has a
presentation
$$u = \prod_{i=1}^e {f_i}^{-1}v^{\pm k}f_i \prod_{i=e+1}^{e+d}{f_i}^{-1}{p_i}^t{p_i}^{-\varphi}f_i,$$
where $f_i \in H*\langle t\rangle, v =
ct\prod_{i=0}^{m}(b_i{a_i}^t)$.

Then in the group $\tilde{G}$ with the starting presentation
$\tilde{G} = \langle G ,t | w^k = 1 \rangle,$ the word $u
=_{\tilde{G}} 1$, $u \in G*\langle t \rangle$ can be represented
as
$$u = \prod_{i=1}^e {f_i}^{-1}w^{\pm k}f_i,$$
indeed, in the starting presentation of $\tilde{G}$ all words of
the type ${f_i}^{-1}{p_i}^t{p_i}^{-\varphi}f_i$ are reduced to the
identity element, and the word $v$ is simply the word $w$, which
is rewritten in new presentation.

Existence of this presentation of $u$ guarantees the existence of
a disk diagram whose exterior face is labelled by $u$ over the
starting presentation of the group $\tilde{G}$. Moreover, let $f'$
be the number of edges of the exterior face of this diagram,
obviously $f' \geq f$, thus $e < \frac{f'}{k-1}$.

Hence we obtain the inequality~(\ref{E:equa4}) with constant $C =
\frac{1}{2(k-1)}$.

\textbf{Case 2.} $|w| =1$. Then the group $\tilde{G}$ has a
presentation:
$$\tilde{G} = \langle G,t | (gt)^k = 1 \rangle \simeq G*\langle x
\rangle_k,$$ Such group is relative hyperbolic with respect to
subgroup $G$.

\section{Proof of corollary \ref{C:SQuniv}.}
Let us show that $G$ is a proper subgroup of the group $\tilde{G}$
and $\tilde{G}$ is non-elementary.

The group $G$ is proper, because the group $\tilde{G}$ has nature
mapping on non-trivial group $\langle t\rangle_k$.

Consider the nature mapping $\tilde{G} \in \hat{G} = \langle
G,t|w=1\rangle$. In \cite{Kl06b} it was proved that the group
$\hat{G}$ contains a nonabelian free subgroup (hence it is
non-elementary), except in the following two cases:
\begin{enumerate}
\item[1)] the group $\hat{G}$ is isomorphic to the Baumslag-Solitar
group $$G_{1,2} = \langle g,t|g^{-1}tg=t^2\rangle;$$
\item[2)] the word $w$ is conjugate to the word $gt$, where $g \in G$.
\end{enumerate}
The group $G_{1,2}$ is non-elementary. So non-elementary quality
of $\tilde{G}$ is followed from non-elementary quality of
$\hat{G}$. In case 2) the group $\tilde{G}$ has a presentation
$G*\langle t \rangle_k$, which is non-elementary.

By the Corollary 1.2 from \cite{AMO07}, we obtain that the group
$\tilde{G}$ is SQ-universal.

\end{document}